# HOMOLOGY OF POLYHEDRA AND QUADRANGULATIONS OF SURFACES

## Serge Lawrencenko

*Department of Mathematics and Science Education*
*Faculty of Service, Russian State University of Tourism and Service*
*Glavnaya 99, Cherkizovo, Moscow Region 141221, Russia*

**Abstract.** A new formula is obtained in algebraic topology, in terms of Betti numbers, and a new method, called the spinal method, is suggested and developed for generating quadrangulations of closed orientable surfaces. Those surfaces arise as the thickenings of 1- and 2-dimensional curvilinear polyhedra, called spines, in Euclidean 3-space. By way of spinal manipulation, quadrangulations with given properties are constructed.

**MSC 2010:** 05C10 (Primary); 57M20, 57M15 (Secondary).

**1. Introduction.** The topics of triangulations and quadrangulations, especially their colorings, have intensively been studied in the literature; cf. [1−6] and the references contained therein. Throughout this paper the term "graph" means "undirected simple finite graph". A *quadrangulation* of the (2-dimensional) sphere with $g$ handles, $\Sigma_g$, with a graph $G$ means a *quadrilateral embedding* $G \to \Sigma_g$ in which each face is a *quadrilateral*—that is, a face bounded by a simple cycle of length 4, and in which each edge meets exactly two quadrilaterals.

In the present paper a class of so-called spinal quadrangulations of $\Sigma_g$ is constructed. The properties of spinal quadrangulations are determined by their spines; especially, the genus $g$ is equal to the 1st Betti number of the spine (Theorem 3.3). This gives rise to the so-called spinal method for generating quadrangulations with prescribed properties. By judiciously choosing the order, the 1st Betti number, and the chromatic number of the spine, one can alter the order, the genus, and the chromatic number of the quadrangulation. For instance, any non-trivial tree is the spine of some bichromatic quadrangulation of the sphere (Corollary 3.5).

New proofs are given for Hartsfield and Ringel's results [7] on minimal quadrangulations. The approach of [7] is based on the method of current graphs which is quite complicated and tedious. In contrast to that method, the spinal method is handy in use and enables us to generate many more minimal quadrangulations (Corollary 3.8). The following is a quadrilateral analogue of Harary, Korzhik and the author's result [1] on triangulations.

**Theorem 1.1.** *Given three integers* $g \geq 0$, $k \geq 2$, $p \geq 4$ *so that:*

(1) $g = 0$ *if* $k = 2$,

(2) $g \geq \frac{1}{2}(k-1)(k-2)$ *if* $k \geq 3$,

(3) $p$ *is even and greater than or equal to* $4g - 2(k^2 - 4k + 2)$,

*there exists a quadrangulation of* $\Sigma_g$ *having $g$ handles, $p$ vertices, vertex chromatic number equal to $k$, and face chromatic number less than or equal to $k$.*

A constructive proof of Theorem 1.1 will be given in Section 4.

**2. Homology of polyhedra.** In this section we motivate, state, and prove a homology lemma in a 2-dimensional setting. Its 1-dimensional variant will be used in Section 3.

Throughout this paper the term "polyhedron" will be used to refer to a *topological (curvilinear) polyhedron*—that is, a topological space of a triangulable set in Euclidean space $E^3$, or in other words, a topological realization $|K|$ of an abstract (finite) simplicial complex $K$ in $E^3$. See Alexandrov [8].

Let $P$ be a compact topological (curvilinear) polyhedron with dimension at most 2 in $E^3$, possibly disconnected. Let $N^3(P)$ denote a 3-dimensional regular neighborhood of $P$ with a small radius. Since $N^3(P)$ is a 3-dimensional solid, its boundary $\partial N^3(P)$ is a possibly disconnected, closed, compact, orientable 2-dimensional surface embedded in $E^3$. Therefore $\partial N^3(P)$ is a disjoint union of a number of topological spheres with handles. Denote the total number of those spheres by $\text{comp}\, \partial N^3(P)$, and the total number of their handles by $\text{hand}\partial N^3(P)$.

The surface $\partial N^3(P)$ can be thought of as the *thickening* of $P$. This operation is used in animation as shown by Brandel, Bechmann, and Bertrand [9]. Therefore construction of quadrangulations on the thickenings of polyhedra may find practical application.

Two surfaces are said to have the same *topological type* if they are homeomorphic. It is not always the case that the topological type of $\partial N^3(|K|)$ is uniquely determined by the original abstract complex $K$. For example, we can add a pair of curved edges to the sphere in a couple of ways: denote by $P_1$ and $P_2$ the polyhedra depicted on the left of Figs. 1 and 2, respectively. In fact $P_1$ and $P_2$ are two non-isotopic topological realizations of the same abstract 2-complex. Then the surface $\partial N^3(P_1)$ is a double torus with a sphere inside it (Fig. 1, right), whereas $\partial N^3(P_2)$ is a torus with another torus inside it (Fig. 2, right).

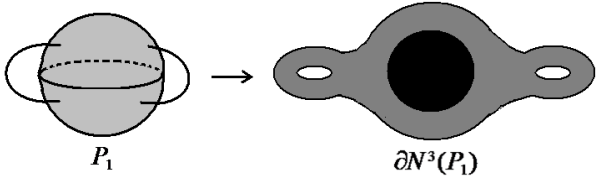

Fig. 1. Double torus with a sphere inside

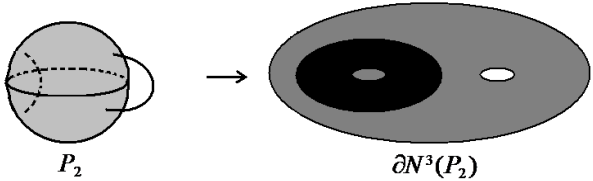

Fig. 2. Torus with another torus inside

Let $\beta_k(P)$ denote the *k*-th Betti number of $P$. Now, we state the homology lemma.

**Lemma 2.1.** *For an arbitrary compact topological (curvilinear) polyhedron $P$ in Euclidean 3-space $E^3$, with dimension 0, 1, or 2, possibly disconnected, the following two identities hold:*

$$\text{comp } \partial N^3(P) = \beta_0(P) + \beta_2(P), \tag{1}$$
$$\text{hand} \partial N^3(P) = \beta_1(P). \tag{2}$$

*Proof.* Let us prove the equalities in the case $P$ is connected. Then it will be easy to generalize to the multi-component situation by summing Eqs. (1), respectively Eqs. (2), over components.

It is a consequence of the Alexander duality theorem that for any compact connected topological polyhedron $P \subset E^3$ with dimension at most 2 the number $\beta_2(P)+1$ is equal to the number of connected (3-dimensional) components of $E^3 - P$. See for instance Zomorodian [10, p. 75]. Then Eq. (1) immediately follows for a connected polyhedron $P$. It remains to prove Eq. (2) for a connected $P$.

If $P$ is a 0-dimensional polyhedron, the surface $\partial N^3(P)$ is a sphere, so Eq. (2) reads $0 = 0$. Eq. (2) also hold in the case $P$ is a 1-dimensional polyhedron—that is, a graph. For let $T$ be a *spanning tree* of $P$—that is, a largest connected subgraph so that $\beta_1(T) = 0$. The surface $\partial N^3(T)$ is homeomorphic to the sphere. It is well-known [11] that $\beta_1(P)$ equals the numbers of edges that are in $P$ but not in $T$. On the other hand, the addition of any such edge to $T$ corresponds to the addition of a handle to $\partial N^3(T)$, and therefore Eq. (2) is still valid.

In the case of a 2-dimensional polyhedron $P$, Eq. (2) is proved by induction on the number of *facets*—that is, topological realizations of the 2-dimensional simplexes of $P$ in $E^3$. The basis of induction corresponds to the situation in which $P$ has no facets at all. Then Eq. (2) holds by the above paragraph.

Assume Eq. (2) holds for any 2-dimensional polyhedron with at most $n$ facets and let $P$ be one with $n+1$ facets. Pick any facet $\sigma^2$ in $P$ and denote by $P - \sigma^2$ the polyhedron obtained by the removal of $\sigma^2$ from $P$, assuming that the boundary $\partial \sigma^2$ remains in $P - \sigma^2$. Eq. (2) holds for $P - \sigma^2$ by the induction assumption.

*Case 1: The 1-dimensional cycle $\partial \sigma^2$ is homologous to zero in $P - \sigma^2$ with respect to integral coefficients.* Hence $\beta_1(P) = \beta_1(P - \sigma^2)$. On the other hand, we have $\beta_0(P) = \beta_0(P - \sigma^2)$ ($= 1$) because both polyhedra are each connected. Furthermore, since the Euler characteristic of $P$ is one greater than that of $P - \sigma^2$, it follows from the Euler-Poincaré formula that $\beta_2(P) = \beta_2(P - \sigma^2) + 1$. It now follows from Eq. (1) (already proved at the beginning of the proof) that

$$\text{comp } \partial N^3(P) = \text{comp } \partial N^3(P - \sigma^2) + 1. \qquad (3)$$

Since $\sigma^2 \cap \partial N^3(\partial \sigma^2)$ is a 1-dimensional circuit (that is, a topological realization of a 1-dimensional cycle in $E^3$) and since $\partial \sigma^2 \subset P - \sigma^2$, it follows that $\sigma^2 \cap \partial N^3(P - \sigma^2)$ is also a 1-dimensional circuit. Therefore $\partial N^3(P)$ results from $\partial N^3(P - \sigma^2)$ by cutting one component of the latter along the circuit $\sigma^2 \cap \partial N^3(P - \sigma^2)$, followed by closing the two resulting holes with 2-disks. Therefore Eq. (3) forces $\text{hand}\partial N^3(P) = \text{hand}\partial N^3(P - \sigma^2)$ after closing the holes.

*Case 2: The 1-dimensional cycle $\partial \sigma^2$ is non-homologous to zero in $P - \sigma^2$ with respect to integral coefficients.* Hence either $\beta_1(P) = \beta_1(P - \sigma^2) - 1$ and we get $\beta_2(P) = \beta_2(P - \sigma^2)$ similarly to case 1, or $\beta_1(P) = \beta_1(P - \sigma^2)$ and we proceed

exactly like in case 1. In the former scenario, by Eq. (1),
comp $\partial N^3(P)$ = comp $\partial N^3(P - \sigma^2)$ and hand $\partial N^3(P)$ = hand $\partial N^3(P - \sigma^2) - 1$. (In fact the latter scenario is impossible in case 2 because the integral homology groups of any finite polyhedron in $E^3$ are known [12, p. 375] to be torsion-free.)

Note that in both cases Eq. (2) remains valid for $P$ and induction completes the proof. ∎

**Remark 2.2.** In fact Lemma 2.1 is a consequence of the following general formula: $\beta_k(\partial N^n) = \beta_k(P) + \beta_{n-k-1}(P)$ ($0 \le k \le n-1$) for any compact topological (curvilinear) polyhedron $P$ with dimension $\le n-1$ in Euclidean $n$-space $E^n$, where $N^n = N^n(P)$ denotes a regular neighborhood of $P$ in $E^n$. This formula is proved by the following standard composition of isomorphisms:

$$\tilde{H}_k(\partial N^n; \mathbf{Q}) \approx \tilde{H}_k(N^n; \mathbf{Q}) \oplus \tilde{H}_k(E^n - N^n; \mathbf{Q}) \approx \tilde{H}_k(P; \mathbf{Q}) \oplus H_{n-k-1}(P; \mathbf{Q}),$$

where the first isomorphism is one from the Mayer-Vietoris exact sequence, and the second is provided by the Alexander duality.

**3. Spinal quadrangulations.** Let $G$ be a graph with $V$ vertices. We associate with $G$ a graph with $2V$ vertices, denoted by $G[:]$ and constructed as follows. Take two disjoint copies, $G'$ and $G''$, of $G$. Vertex $v'$ in $G'$ and the corresponding vertex $v''$ in $G''$ are called the *twins* of each other. Every vertex in $G'$ is joined by edges to all neighbors of its twin vertex in $G''$, but not to the twin itself. More precisely, for each vertex $v$ of $G$ the vertices $v'$ and $v''$ are non-adjacent in $G[:]$, and for each edge $vu$ of $G$ the graph $G[:]$ has four edges as follows: $v'u'$, $v'u''$, $v''u'$, and $v''u''$. In fact, $G[:] = G[\overline{K_2}]$, where the inner factor $\overline{K_2}$ in the composition (see [11] or [14] for the definition) is the graph complementary to the complete graph $K_2$, or in other words, $\overline{K_2}$ is the graph consisting of two isolated vertices (which is what the colon in the notation $G[:]$ is meant to indicate). We call the graph $G[:]$ the *2-fold interlacement* of $G$. Observe that the number of edges in $G[:]$ is four times the number of edges in $G$.

**Lemma 3.1.** *For any graph $G$ without isolated vertices and possibly disconnected, its 2-fold interlacement quadrangulates its thickening. In other words, $G[:]$ can be embedded on $\partial N^3(G)$ with each face a quadrilateral.*

***Proof.*** It suffices to prove the statement for any one connected component of $G$, so we may assume that $G$ is connected. Prepare a set $\{S_v\}$ of disjoint spheres respectively corresponding to the vertices $v$ of $G$. For each $v$ mark two antipodal points on $S_v$ to represent its twins $v'$ and $v''$, and cut out $\deg v$ quadrilaterals from $S_v$, where $\deg v$ denotes the *degree* of $v$ in $G$ —that is, the number of edges that meet at $v$. The resulting quadrilateral holes are shaded black in Fig. 3. We then repeatedly loop through all the edges in $G$ and, for each current edge $vu$, arbitrarily choose a pair of holes (not yet chosen before), one in $S_v$, denote it by $h_{vu}$, and the other one in $S_u$, denote it by $h_{uv}$.

For each $v$ and for each hole in $S_v$, draw four arcs on $S_v$ as fragments of the edges of $G[:]$ to be constructed, as shown in Fig. 3. Two of those arcs join $v'$ to the upper two corners of the hole (those corners are not vertices!), two more join $v''$ to the lower two corners.

It remains to close the holes to get a desired quadrangulation $G[:] \to \partial N^3(G)$. We now show, for each edge $vu$ of $G$, how to close a pair of holes $h_{vu}$ and $h_{uv}$. Prepare a topological parallelepiped tube with open ends, with its four long edges representing the four remaining fragments of the edges under construction. Now, to complete the edges, add the tube to connect the boundaries of $h_{vu}$ and $h_{uv}$, giving a quarter-twist to one end of the tube in either direction (either clockwise or counterclockwise). For certainty, the direction is chosen counterclockwise as indicated by the curved arrow in Fig. 3.

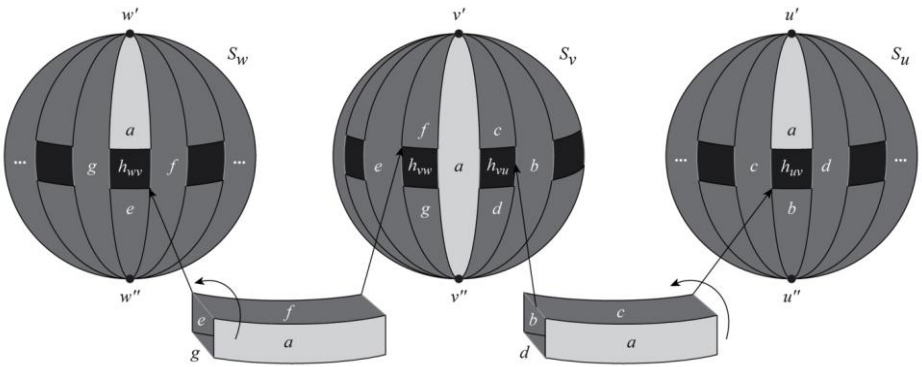

Fig. 3. Diagram for the Proof of Lemma 3.1

To keep a track of the faces appearing in the resulting embedding, the areas sharing common border with $h_{vu}$ and $h_{uv}$ are marked with letters $a$, $b$, $c$, and $d$ in Fig. 3. Face $a$ is shaded lighter than the rest of the surface for convenience of observation. In the case $\deg v = 1$ there is only one hole in $S_v$, say $h_{vu}$, in which case faces $a$ and $b$ merge into one. We close the holes $h_{vu}$ and $h_{uv}$ by adjoining the tube onto their boundaries. The quarter-twist of the tube between $S_v$ and $S_u$ ensures that the areas with the same letter merge to form one face. The marking on $S_u$ in Fig. 3 corresponds to the choice of the counterclockwise quarter-twist in the preceding paragraph. If we had chosen the clockwise quarter-twist, we would have switched around the letters $a$ and $b$ on $S_u$ as well as $c$ and $d$.

Clearly, with all the holes closed, we obtain a surface homeomorphic to $\partial N^3(G)$. The graph obtained on the surface is indeed the 2-fold interlacement of $G$ because for each vertex $v$ of $G$ the vertices $v'$ and $v''$ are non-adjacent, and for each edge $vu$ of $G$ the graph constructed on the surface has four edges as follows: $v'u'$, $v'u''$, $v''u'$, and $v''u''$ (check with Fig. 3).

It remains to prove that the embedding $G[:] \to \partial N^3(G)$ constructed is indeed a quadrangulation. We show it in two ways: (I) directly, and (II) via the Euler equation along with Lemma 2.1.

(I) Here we identify a typical face of the embedding. For certainty, consider the area marked by $a$, shaded light gray in Fig. 3. Let $h_{vw}$ denote the other hole sharing a common border with area $a$ on $S_v$. Then $a$ extends onto $S_w$ similarly to how it extends onto $S_u$. See Fig. 3. It is not hard to traverse the boundary of $a$ and find that that boundary is in fact a cycle of length 4: $(v', u', v'', w')$. A general face may have one of the following four forms: $(v', \hat{u}, v'', \hat{w})$, where $\hat{u} \in \{u', u''\}$ and $\hat{w} \in \{w', w''\}$, depending on the two possible directions (that is, clockwise or counterclockwise) of the quarter-twists given to the tube between $S_v$ and $S_u$, and to the tube between $S_v$ and $S_w$ (respectively). Anyhow, each face is indeed a quadrilateral.

(II) Now we show that the embedding $G[:] \to \partial N^3(G)$ is a quarangulation by the Euler equation. Let $F$ denote the number of faces in this embedding and let $V$ and $E$ denote the numbers of vertices and edges in $G$, respectively. Thus the 2-fold interlacement $G[:]$ has $2V$ vertices and $4E$ edges. Then the Euler equation for $G[:] \to \partial N^3(G)$ reads as

$$2V - 4E + F = 2 - 2g, \qquad (4)$$

where $g$ is the total number of handles in $\partial N^3(G)$. Therefore, thanks to Lemma 2.1,

$$g = \text{hand} \partial N^3(G) = \beta_1(G) = E - V + 1, \qquad (5)$$

where the last equation in Eqs. (5) is a well-known expression for the first Betti number of a connected graph (see [11, 14]). We derive from Eqs. (4) and (5) that $2 \times (4E) = 4 \times F$ —that is, twice the number of sides (edges) is equal to four times the number of faces in the constructed embedding $G[:] \to \partial N^3(G)$. Hence each face is indeed a quadrilateral since the same edge is always a side in two faces on any closed surface. The proof is complete. ∎

The *genus of* an arbitrary (possibly disconnected) *graph* $G$ is the minimum $g$ so that there exists an embedding $G \to \Sigma_g$; such an embedding is called a *minimum-genus embedding*. In the case when $G$ has no cycle of length 3, the graph $G[:]$ has no such cycle either. Hence any embedding of $G[:]$ has no triangular face. Hence any quadrangulation with the graph $G[:]$ would be automatically a minimum-genus embedding (see White [14, p. 54]). Therefore, since $G[:] = G[\overline{K_2}]$, a classical result by White [13, 14] follows, as a byproduct, from a combination of Lemma 3.1 and Eq. (2):

**Corollary 3.2.** *For any (possibly disconnected) graph $G$ without isolated vertices and without cycles of length 3, the genus of the graph $G[\overline{K_2}]$ is equal to $\beta_1(G)$.*

From now until the end of the paper we assume that $G$ is a connected graph and therefore $\partial N^3(G) = \Sigma_g$ for some $g$. A *spinal quadrangulation* is denoted by $Q = Q(G[:])$ and is defined to be any quadrilateral embeddings $Q : G[:] \to \Sigma_g$ that can be obtained by the construction in the proof of Lemma 3.1. The graph $G$ is called the *spine of $Q$*. The number of handles $g$ is called the *genus of $Q$*. It should be noted that the construction procedure retains some freedom of action (since the holes can be grouped into pairs in different ways) and therefore there are, in general, many spinal quadrangulations with the same spine.

Lemmas 2.1 and 3.1 imply the following crucial result.

**Theorem 3.3.** *The genus of a spinal quadrangulation is equal to the 1st Betti number of the spine.* ∎

A *proper vertex coloring* of a graph $G$ is defined to be a coloring of the vertices of $G$ so that no two adjacent vertices share the same color. Similarly, a *proper face coloring* of a quadrangulation $Q$ is a coloring of the quadrilaterals of $Q$ so that any two adjacent quadrilaterals are colored differently, where two quadrilaterals are said to be adjacent if they meet the same edge. The minimum number of colors sufficient for a proper vertex coloring of $G$ is called the *vertex chromatic number* of $G$ and is denoted by $\chi_V(G)$. Similarly, the least number of colors in a proper face coloring of $Q$ is called the *face chromatic number* of $Q$ and is denoted by $\chi_F(Q)$.

**Lemma 3.4.** $\chi_V(G[:]) = \chi_V(G) \geq \chi_F(Q(G[:]))$.

**Proof.** Firstly, we prove that $\chi_V(G[:]) = \chi_V(G)$. In fact, to get a proper vertex coloring of $G[:]$ with exactly $\chi_V(G)$ colors, fix any proper vertex coloring of $G$ with this number of colors, and then color the vertices of $G[:]$ by applying the color of $v$ to both twins $v'$ and $v''$ in $G[:]$, for each vertex $v$ of $G$. We thus indeed obtain a proper vertex coloring of $G[:]$, which immediately follows from the definition of $G[:]$.

On the other hand, $G[:]$ cannot be vertex colored with less than $\chi_V(G)$ colors because $G[:]$ contains $G$ as a subgraph.

We proceed to prove the inequality stated in the lemma. Let us say that a vertex $v$ of the spine $G$ is the *source* for a quadrilateral of $Q(G[:])$ if that quadrilateral meets both twins $v'$ and $v''$. It is clear that each quadrilateral has one and only one source. For instance, the typical quadrilateral $a$ has $v$ as its source (check with Fig. 3). This quadrilateral is adjacent to quadrilaterals $c$ and $d$ with source $u$, and to quadrilaterals $f$ and $g$ with source $w$. Therefore, any two quadrilaterals with the same source are non-adjacent, and two quadrilaterals are adjacent in $Q(G[:])$ only if their sources are adjacent as vertices of $G$. It follows that, given any proper vertex coloring of the spine $G$, we obtain a proper face coloring of $Q(G[:])$ by applying to each quadrilateral the color of its source. ∎

The following is a corollary which may be used in applications of graph theory.

**Corollary 3.5.** *For every non-trivial tree $T$, the graph $T[:]$ quadrangulates the sphere with both vertex and face chromatic numbers equal to 2.* ∎

**Corollary 3.6.** *For every integer $n \geq 2$ the 2-fold interlacement $K_n[:]$ of the complete graph $K_n$ quadrangulates $\Sigma_g$ with $g = \frac{1}{2}(n-1)(n-2)$ and vertex chromatic number equal to $n$, and face chromatic number less than or equal to $n$.*

*Proof.* The spinal method works as follows. Compute

$$\beta_1(K_n) = \# \text{ of edges} - \# \text{ of vertices} + 1 = \tfrac{1}{2}n(n-1) - n + 1 = \tfrac{1}{2}(n-1)(n-2),$$

and apply Theorem 3.3 and Lemma 3.4. The proof is complete. ∎

**Corollary 3.7.** *Let $n$ be an integer $\geq 3$ and let $K_n - e$ be the graph obtained from $K_n$ by removing one edge (without incident vertices). Then $(K_n - e)[:]$ quadrangulates $\Sigma_g$ with $g = \tfrac{1}{2}\left[(n-1)(n-2) - 2\right]$ and vertex chromatic number equal to $n-1$, and face chromatic number less than or equal to $n-1$.*

*Proof.* Like in the proof of Corollary 3.6, compute

$$\beta_1(K_n - e) = \left[\tfrac{1}{2}n(n-1) - 1\right] - n + 1 = \tfrac{1}{2}\left[(n-1)(n-2) - 2\right],$$

and the statements now follow from Theorem 3.3 and Lemma 3.4. ∎

The quadrangulations of Corollaries 3.6 and 3.7 were first discovered by Hartsfield and Ringel [7] and, when $n$ is sufficiently large, were shown to be minimal quadrangulations. A quadrangulation of a fixed surface $\Sigma_g$ is said to be *minimal on* $\Sigma_g$ if the number of vertices is minimal among all quadrangulations of $\Sigma_g$. The spinal method enables us to go further and generalize the results of [7] by way of generating many more minimal quadrangulations as follows.

**Corollary 3.8.** *Let $n$ and $m$ be positive integers so that $n \geq 2$ and $m \leq n-1$, and let $K_n - E(K_m)$ be the graph obtained from $K_n$ by removing the edges (without incident vertices) of some subgraph $K_m$. Then the 2-fold interlacement of $K_n - E(K_m)$ quadrangulates $\Sigma_g$ with $g = \frac{1}{2}[(n-1)(n-2) - m(m-1)]$ and vertex chromatic number equal to $n - m + 1$, and face chromatic number less than or equal to $n - m + 1$. Furthermore, for $g \geq 1$, any such quadrangulation is minimal on $\Sigma_g$ whenever $n \geq 4 + 2m(m-1)$.*

*Proof.* This proceeds by the spinal method used in the proofs of Corollaries 3.6 and 3.7. Compute

$$\beta_1(K_n - E(K_m)) = \left[\tfrac{1}{2}n(n-1) - \tfrac{1}{2}m(m-1)\right] - n + 1 = \tfrac{1}{2}[(n-1)(n-2) - m(m-1)],$$

and the statement follows.

It remains to prove minimality of the spinal quadrangulation constructed. Let $V$, $E$, $F$ denote the number of vertices, edges, and quadrilaterals of an arbitrary quadrangulation of $\Sigma_g$. Then the Euler equation reads as

$$V - E + F = 2 - 2g.$$

Furthermore, since any pair of vertices are joined by at most one edge, we have $E \leq \binom{V}{2}$, and since each edge meets exactly two quadrilaterals, we have $4F = 2E$. From these it can be derived that $V^2 - 5V + (8 - 8g) \geq 0$. This quadric inequality has the solution:

$$V \geq \left\lceil \tfrac{1}{2}\left(5 + \sqrt{32g - 7}\right)\right\rceil. \tag{6}$$

(It should be noted that this inequality is valid if and only if $g \geq 1$.) In the beginning of the proof we have proved that for the spinal quadrangulation constructed with the spine $K_n - E(K_m)$ we have $g = \tfrac{1}{2}[(n-1)(n-2) - m(m-1)]$, hence

$$32g - 7 = 16n^2 - 48n + 25 - 16m^2 + 16m. \tag{7}$$

Now it follows from Eqs. (6) and (7) that the constructed quadrangulation is minimal on $\Sigma_g$ whenever the following double inequality holds:

$$2n - 1 < \tfrac{1}{2}\left(5 + \sqrt{16n^2 - 48n + 25 - 16m^2 + 16m}\right) \leq 2n.$$

This can be rewritten as

$$4n-7 < \sqrt{16n^2 - 48n + 25 - 16m^2 + 16m} \le 4n - 5,$$

or as

$$16n^2 - 56n + 49 < 16n^2 - 48n + 25 - 16m^2 + 16m \le 16n^2 - 40n + 25,$$

or as

$$\begin{cases} 8n > 24 + 16m(m-1), \\ 8n \ge -16m(m-1). \end{cases}$$

Since the second inequality is a tautology, we conclude that $n > 3 + 2m(m-1)$, hence $n \ge 4 + 2m(m-1)$. The proof is complete. ∎

**4. Proof of Theorem 1.1.** By Theorem 3.3, to satisfy the requirement on the number of handles, it suffices to take as a spine any graph $G$ with $\beta_1(G) = g$. To construct a spine to satisfy all the requirements for $g = 0$, we take as a spine any tree with $p/2$ vertices and we are done by Corollary 3.5. For $g \ge 1$, we start up with the complete graph $K_k$ with $k \ge 3$. Hypothesis (2) of the theorem guarantees that $\beta_1(K_k) \le g$. If $\beta_1(K_k) = g$, the quadrangulation $Q(K_k[:])$ already has the required number of handles and the required chromatic numbers by Corollary 3.6. If $\beta_1(K_k) < g$, we hang up $g - \beta_1(K_k)$ edges (as shown in Fig. 4) until we reach the total of $g$ independent cycles required. By now, the number of vertices in the spine under construction equals to $k + 2[g - \beta_1(K_k)] = 2g - (k^2 - 4k + 2) \le p/2$, where the upper bound is provided by hypothesis (3). Finally, if this bound has not been attained, we hang up necessary vertices until we reach the total of $p/2$ vertices required in the spine. The construction is illustrated in Fig. 4, where a possible spine is shown for the triple $g = 5$, $k = 4$, $p = 20$. The proof is complete.

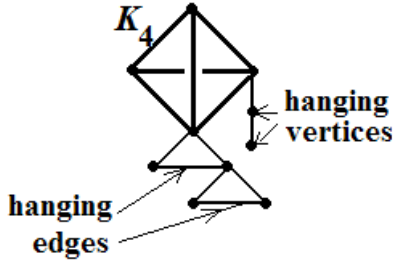

Fig. 4. Diagram for the Proof of Theorem 1.1

**Acknowledgments.** I thank A.M. Gurin, V.Yu. Rovenski, I.Kh. Sabitov, and N.P. Strelkova for productive discussions.